\documentstyle[12pt,fleqn]{article}
\renewcommand{\baselinestretch}{1.3}
\newcommand{\be}{\begin{equation}}
\newcommand{\ee}{\end{equation}}
\newcommand{\bea}{\begin{eqnarray}}
\newcommand{\eea}{\end{eqnarray}}
\newcommand{\su}{\mbox{$SU_q(2)$\,\,}}
\newcommand{\pro}{\partial}
\newcommand{\op}{\omega^{++}}
\newcommand{\om}{\omega^{--}}
\newcommand{\oo}{{\omega^0}}
\newcommand{\dfrac}{\displaystyle\frac}
\newcommand{\ba}{\begin{array}}
\newcommand{\ea}{\end{array}}
\newcommand{\dpst}{\displaystyle}
\newcommand{\dpp}{{\cal D}^{++}}
\newcommand{\dmm}{{\cal D}^{--}}
\newcommand{\docal}{{\cal D}^{0}}

\newcommand{\vareps}{\varepsilon}

\newcommand{\qin}{\displaystyle\frac{1}{q}}
\newcommand{\uone}{\mbox{$U(1)$\,\,}}
\newcommand{\suu}{\mbox{$SU(2)$\,\,}}
\newcommand{\nn}{\nonumber}
\newcommand{\pbar}{\bar{\phi}}
\begin{document}
\large
\title{  DIFFERENTIAL CALCULI ON THE QUANTUM GROUP
             \su AND GLOBAL  $U(1)$-COVARIANCE}

\author{ D.G. Pak \thanks {E-mail: dmipak@silk.glas.apc.org}}
\date{28 August 1995}
\maketitle
\begin{center}
       Department of Theoretical Physics, Research Institute\\
       of Applied Physics, Tashkent State  University,\\
       Vuzgorodok, 700095, Tashkent, Republic of Uzbekistan \\
\end{center}
\begin{abstract}
      A variety of three-dimensional left-covariant differential
calculi on the quantum group \su is considered using an approach
based on global $U(1)$-covariance. Explicit representations of
possible $q$-Lie algebras are constructed in terms of
differential operators. A gauge covariant
differential algebra is uniquely determined. The non-standard
Leibnitz rule is obtained for a corresponding $q$-Lie algebra.
\end{abstract}
   \vskip 50pt

\newpage
\section{Introduction}
\indent

     One of the features of non-commutative geometry in the
quantum group theory [1-5]
 is non-uniqueness
in defining a differential calculus on the quantum groups and quantum spaces.
The bicovariance condition determines a unique differential calculus on the
linear quantum groups $GL_{q}(N)\,$ (up to symmetry corresponding to the
exchange $q \rightarrow \qin$) [6, 7]
 and
provides existence of the corresponding
 gauge covariant differential algebra \cite{Isa2}.
Direct reducing the $GL_{q}(N)$-bicovariant differential calculus to
a case of the special linear quantum group $SL_{q}(N)\,$
encounters  difficulties connected with a loss
 of the centrality condition for
a quantum determinant. Four-dimensional $4D_{\pm}$ \, bicovariant and
three-dimensional $(3D)$ left-covariant differential calculi on the
simplest special unitary
quantum group \su were considered as well using a
standard Woronowicz approach [3, 6].
 A full consistent construction of the 3D bicovariant differential calculus
 and a gauge covariant differential algebra on the \su
are unknown up to now, furthermore, there are strong
limitations imposed by no-go theorems \cite{Are1}.
 A possible way to solve this
 problem suggests using a non-standard Leibnitz rule as it was
pointed out in ref. \cite{Fad2}.

    In this paper possible 3D left-covariant differential calculi and
gauge covariant differential algebra on the quantum group \su
are considered in the framework of approach which respects a global
 $U(1)$-covariance.
 The group \uone is a stabilizing subgroup for the quantum
 group \su and the $U(1)$-covariant treatment allows to pass
straightforwardly to the description of the quantum sphere
 $S_{q}^{2} \, \sim \,\newline  {\suu } / {\uone}$. In Section 2 we
 construct
explicit representations of $q$-Lie algebras of left-invariant
vector fields on \su in terms of differential operators.
The $U(1)$-covariance constraint reduces the variety of possible
left-covariant
differential calculi on \su and leads to a unique gauge covariant
differential algebra as it is shown in Section 3. The main
commutation relations for the differential algebra agree
with ones obtained earlier in the bicovariant
formalism [8, 10].
The principal difference of our approach is that we propose a
non-standard Leibnitz rule which is consistent with the
gauge $SU_q(2)$-covariance.
 Some discussion of a quantum group gauge Yang-Mills
theory is given in Conclusion.

\section{Left-invariant vector fields}
\indent

  Following the $R$-matrix formalism \cite{Fad1} the main commutation
relation for the generators $T^i_j \;\;
(i, j =1,2)$ of the quantum group \su are defined by a
standard $R$-matrix as follows
\be
   R_{12}T_1T_2 = T_2T_1R_{12}.
\ee
Let us choose a covariant parametrization for the matrix
 $T^i_j$
\bea
 T^i_j=\left(\begin{array}{cc}
                  y^1&x^1                              \label {param}    \\
                  y^2&x^2
                 \end{array} \right)
                             \equiv (y^ix^i),
\eea
where  $x^i, y_i \  $ are generators (coordinates)
of the function algebra on the quantum hermitean
vector space $U^2_q$  endowed with an involution
${\ast}:{ \stackrel{\ast}{x^{i}} = y_i}$
and $SU_q(2)$-comodule structure.The unimodularity condition
 takes a simple covariant form
\bea
\ba{lr}
  {\cal D} \equiv {\det}_{q} T^i_j = x_i y^i = 1,&
                     x_i=\vareps_{ij}x^j.
\ea
\eea
 The \su indices are raised and lowered with the
invariant metric $ \varepsilon_{ij} \,\,(\varepsilon_{12}=1,  \,\,
\varepsilon_{21}=-q^{-1})$.

     The parametrization (\ref{param}) was used in a harmonic
formalism \cite{Gal1} of extended superfield supesymmetric theories. The
coordinates $(x, y)$ parametrize the
quantum sphere
 $S^{2}_q  \, \sim \,  {\suu } / {\uone}$
and are just the quantum generalizations of classical
 harmonic functions $(u^{\pm})$ (so called "harmonics")
 \bea
     x^i \equiv u^{+i}, \,\,\,\,\,\,\,y^i\equiv u^{-i}.
\eea
The signs $(\pm)$ correspond to  charges $(\pm1)$
of the stabilizing subgroup \uone for the quantum
group \su. To simplify notations we shall not pass to the
notations adopted in the harmonic formalism keeping in mind
that all geometric objects (like coordinates, derivatives, differential
forms
 etc.) have definite \uone charges.

    Consider main commutation relations between the coordinates  \newline
$(x,y)$ and derivatives
 $\partial_i\equiv \dfrac{\partial}{\partial
 x^i}, \;\;
 \bar\partial^i\equiv \dfrac{\partial}{\partial y_i}$ on the
quantum group \su:
 \bea
\ba{lr}
R_{12}(\pro_T)_1 (\pro_T)_2 = (\pro_T)_2 (\pro_T)_1 R_{21} ,&
    (\pro_T)^i_j \equiv \left(\begin{array}{cc}
                  \bar \pro_1 &\bar \pro_2\\
                  \pro_1 &\pro_2
                 \end{array} \right)     ,
\ea
\eea
\bea
\ba{lr}
    \pro_ix^k=\gamma_i^k+qY^{nk}_{mi}x^m\pro_n,  &
                        {\bar\pro}^i\bar y_j=\delta^i_j+qy_m\bar\pro^n\hat
                        R^{mi}_{nj},  \\             \label{pro}
    \pro_i y_j=q(\hat R^{-1})^{lk}_{ji} y_k\pro_l,   &
                        \bar \pro^i x^j = \qin \hat R^{ij}_{kl}
                        x^k\bar\pro^l,
\ea                                        \label {comm}
\eea
here, we use standard definitions for the matrices $\hat R^{ij}_{kl},\;
    Y^{ij}_{kl}\;$.
The commutation relations (\ref {comm}) do
not differ on principle from
ones given in ref. \cite {WZ1}. Our choice is
motivated by applying manifest covariant tensor notations which
 are convenient
in constructing explicit representations for the $q$-Lie algebras.
Thus, one implies all geometric objects with upper (lower)
indices to be transformed under the quantum group co-action
$ \Delta $ like classical co- (contra-) variant tensors.
For instance, a second rank tensor $N_i^j$
will be transformed as follows
\be
      (N_i^j)^{\prime} = (T^{\dag})^k_i {T^j}_l N_k^l
\ee
(Hereafter the signs $ \bigotimes$ of tensor product are omitted).

      Let us define the left-invariant first-order differential operators
\bea
    D^{++}\equiv x_i\bar\pro^i,   &   D^{--}\equiv -y_i\pro^i,
\eea
where ($\pm \pm$) correspond  to \uone charges ($\pm 2$).
The action of the operators
  $ D^{\pm\pm}$ on the coordinates $(x, y)$ has a simple form
\bea
    D^{++}x^i=0,  &    D^{--}x^i=y^i,       \label {def1}    \nonumber\\
    D^{++}y_i=x_i,  &  D^{--}y_i=0.
\eea
 The Leibnitz rule for these differential operators
 may be written in a  convenient form if one
considers their action  on functions with definite
\uone charges. After some calculations one finds
\be
   D^{\pm\pm}(f^{(m)} g^{(n)}) = (D^{\pm\pm}f^{(m)}) g^{(n)}+
                         q^{-m}f^{(m)} D^{\pm\pm}g^{(n)}.
                                             \label {def2}
\ee

     This is a special feature of quantum group non-commutative geometry
that the quantum analogue to  classical
\uone generator
 can be realized as a second-order differential operator.
With a little algebra one can write the next expression
for the quantum \uone generator
\be
   D^0\equiv -x_i\pro^i-q^2 y_i\bar\pro^i +
                    (1-q^2)x_iy_k\bar\pro^k\pro^i.
\ee
The operator  $D^0$ has eigenfunctions which are just the functions with
definite \ \uone \  charges
\bea
   D^0f^{(n)} = {\{n\}}_q f^{(n)},
          \label {def3}\\
 \{ n\} _q \equiv
   \dfrac {1-q^{-2n}}
                          {1-q^{-2}} ,                            \nonumber
\eea
where    $ {\{n\}}_q $ \  is a $q$-number. It is not hard to check the
following Leibnitz rule for the operator
 $D^0$
\be
   D^0(f^{(m)}g^{(n)}) = (D^0f^{(m)})g^{(n)} + q^{-2m}f^{(m)}
                                     D^0g^{(n)}.           \label {def4}
\ee
Reducing the space of functions on \su to the space of functions
 with a definite  \uone charge one obtains
 the covariant description of the
coset $S_{q}^{2} \, \sim \,  {\suu } / {\uone}$.

By direct calculating one can verify that the operators
$D^{\pm\pm, 0}$ form the generalized $q$-Lie
algebra of \su \ \cite {Schi1}
\bea
\ba{lr}
  {[D^0,D^{++}]}_{\dpst {q^{-4}}} = {\{ 2\}} _q D^{++}, &  \\
   {[D^0, D^{--}]}_{\dpst {q^4}} = {\{ -2\}} _q D^{--},  & \label{Dalg}   \\
   {[D^{++}, D^{--}]}_{\dpst {q^2}} = D^0,                     &
\ea
\eea
here, $[A,B]_{\dpst {q^s}} \equiv AB-q^s BA$.
Note, that the algebra (\ref{Dalg}) is valid
irrespective of whether one imposes the unimodularity
constraint ${\cal D}=1$.
Observe that the action of the operators
$D^{\pm\pm, 0}$  is consistent with the constraint
${\cal{D}} = 1 $ , so  we have
\be
   D^{\pm\pm0}({\cal D} - 1 )\: f(x,y))\cong 0.
\ee
The symbol $\cong$ means that one has a weak equality which is
fulfilled in virtue of commutation relations.
We shall treat the algebra (\ref{Dalg}) as a main $q$-Lie algebra
of left-invariant vector fields on the quantum group \su.
A corresponding $q$-generalized Jacobi identity is available
\bea
    [D^0, [D^{++}, D^{--}]_{\dpst {q^2}}]
 + [D^{++},[D^{--},D^0]_{\dpst {q^{-4}}}]_{\dpst {q^{-2}}}   \nn  \\
            +q^2[D^{--},[D^0,D{++}]_{\dpst {q^{-4}}}]_{\dpst {q^{-2}}} = 0.
\eea

Let us  now  pass to constructing other possible $q$-Lie algebras of
left-invariant vector fields on the \su. For this purpose we consider
 differential operators $\mu,\,\,\nu\,\,$ \cite{Ogi1}
 \be
   \mu = 1+(q^2-1) y_i{\bar {\pro}}^i,  \,\,\,\,\,\,\,
   \nu = 1+ (1- \frac{1}{q^2}) x_i \pro^i.
\ee
One can then see that the operators
 $\mu,\,\,\nu\,\,$ obey the simple commutation
relations
\bea
   \mu D^{--} = q^2 D^{--} \mu,  & \,\,\,\,\, \nn
   \mu D^{++} = \dfrac {1}{q^2} D^{++}\mu, \\
   \mu D^0 = D^0 \mu,            & \,\,\,\,\,  \mu\nu=\nu\mu  .
\eea
Similar formulae hold for the operator $\nu \,$ as well.
Using these relations one can find that
 the operators
    $\dpp,\,  \dmm ,\, \docal $ defined by the next equations
\bea
\ba{lr}
    \dpp = \mu^{-\frac {1}{2}} D^{++},& \dmm = \nu^{-\frac{1}{2}} D^{--} ,\\
    \docal = \qin \mu \nu D^0  \equiv {[\pro^0]}_q ,   &
\ea
\eea
 generate just the Drinfeld-Jimbo quantum enveloping algebra.
 To construct  other possible $q$-Lie algebras
one introduces another  differential operators
$\Delta^{++}, \,\Delta^{--},\, \Delta^{0}  \,$ as follows
\bea
\ba{lr}
\Delta^{++}=F(\hat Z)D^{++}, &  \\
\Delta^{--}=G(\hat Z)D^{--},& \\
\Delta^{0}=H(\hat Z), & \\
\hat Z \equiv (\mu \nu)^{-\frac{1}{2}}, & \hat Zf^{(n)}=q^n f^{(n)} ,
\ea
\eea
here,  $F,\, G,\, H$  --
are some operator  functions of $\hat Z$. The operators $\Delta^{\pm\pm,0}$
generate a quantum enveloping
 algebra of left-invariant vector fields
with  $q$-generalized commutators:
\bea
\Delta^{++}\Delta^{--} - q^{2p}\Delta^{--}\Delta^{++} =
\dfrac{q^p}{q-q^{-1}} {\hat Z}^{-p}(\hat Z-{\hat Z}^{-1}),  \label {quad}\\
\Delta^{0}\Delta^{++}-q^{2s}\Delta^{++}\Delta^{0}=\Delta^{++},
\,\,\,\,\,\,\,\,\,\,\,\,\,\,\,\,\,\,
\,\,\,\,\,\,\,\,\,\,\,\,\,\,\,\,\,\,     \\
\Delta^{0}=\dfrac{1-{\hat Z}^s}{1-q^{2s}},  \label{psalg}
\,\,\,\,\,\,\,\,\,\,\,\,\,\,\,\,\,\,\,\,\,\,\,\,\,\,\,\,\,\,\,\,
\,\,\,\,\,\,\,\,\,\,\,\,\,\,\,\,\,\,\,\,\,\,\,\,\,\,\,\,\,\,\,\,
\eea
where $s,p$
    --  are arbitrary integers which determine in part
the functions $F, G$. Using the last equation (\ref{psalg}) we may express
the operator $\hat Z$ \  in terms of $\Delta^{0}$  and then restrict
 the arbitrariness of the parameters $s, p$
 by considering only quadratic relations.
 Having put the quadratic terms from
 the r.h.s. of the eqn. (\ref {quad}) to the l.h.s. one obtains the
possible 3D $q$-Lie algebras of \su.

    A special choice of a $q$-Lie algebra with some assigned Leibnitz rule
 defines uniquely the
differential calculus on the quantum group \su. Exterior differential
1-forms are treated as dual elements to the left-invariant
 vector fields.

\section{Left-covariant differential algebras}
\indent

    In this section we give description of possible  \su left-covariant
differential algebras with \uone conserved charge. The gauge covariance
condition leads to a unique differential algebra of \su. At the same time
a Leibnitz rule for the exterior differential is not fixed yet.
To find the differentiation rules one needs to choose a corresponding
 $q$-Lie algebra of left-invariant vector fields.

   Let us consider the left-invariant Cartan 1-forms
$\Omega$ on the quantum group \su
\bea
\Omega=dT^{-1} T \equiv \left(\begin{array}{c cr}
        \oo &\op  \\
        \om &-q^{2}\oo
                       \end{array}
                                  \right),
\eea
where $ \oo, \op, \om $ are the basic left-invariant differential 1-forms
with corresponding  \uone charges $( 0, +2, -2)$. One defines
a gauge transformation as follows
\bea
    T^{g}=\tilde {T} T, \,\,\,\,\,\,\,\,\,\,\,\,\,\,\,\,\,\,\,\,\,\,\nn  \\
   \Omega^{g} = \Omega - T^{-1} \tilde {\Omega}T,  \,\,\,\,\,\,\,\,
                                          \tilde{\Omega} \equiv d\tilde
                                        {T}^{-1}\tilde{T}.
\eea

It turns out that the requirement of the
 global $U(1)$-covariance and
the consistence with the quantum group structure determine
uniquely all commutation relations
between the differential 1-forms $\omega$ and the coordinates $(x,y)$.
As a result we have
\bea
\ba{ll}
\op x=qx\op , &   \om x= \qin x\om+ \dfrac{1-q^{4}}{q} y\oo,     \\
\op y=\qin y\op, &  \om y=q y \om,\\
\oo x= x\oo + (1-\dfrac{1}{q^2}) y \op, &     \\
\oo y=y\oo .   &
\ea
\eea

Similar consideration of commutation relations for the basic
 differential 1-forms ${\omega}^{\pm\pm,0}$ leads to the left-covariant
algebras parametrized by a real number $\sigma$:
\bea
\op \op =\om \om = 0,   \label {f4alg1}  \\
\op \oo + q^2 \oo \op =0, \label{f4alg2}  \\
\om \oo + \dfrac{1}{q^2} \oo \om =0,     \label {f4alg3} \\
\op \om + q^{\sigma} \om \op +\dfrac {q^2 (1-q^{\sigma})(1+q^2)}{q^2-1}
                       \oo \oo =0,  \label{f4alg4}\\
\oo \oo = \dfrac{1-q^2}{q^2 (1+q^2)} \op \om.       \label {f4alg5}
\eea

It should be noted that the algebra defined by eqs.
(\ref{f4alg1}-\ref{f4alg5})
 is left-covariant irrespective of  whether one considers
  the last relation (\ref {f4alg5}).
  Requiring the covariance under the gauge transformations and
  using the additional
  commutation relation
  \bea
\tilde {\Omega} \Omega =-q^2 \Omega \tilde {\Omega}
\eea
   one finds a unique gauge covariant differential algebra
at   $\sigma =4$
  \bea
\ba{lr}
\op \op =\om \om = 0, &\\
 \op \oo + q^2 \oo \op =0,&   \\
 \om \oo + \dfrac{1}{q^2}
                    \oo \om =0, &    \label{qalg}       \\
(1+q^2)^2 \oo\oo =
 \dfrac{1}{q^2}\op\om+q^2 \om\op.    &
\ea
\eea

The equation (\ref {f4alg5}) is not gauge covariant and should be omitted.
 So defined gauge covariant differential
algebra coincides with one obtained in  ref.
\cite {Isa2}, where possible $GL_q(N)$-covariant quantum algebras
were studied. It should be noted that
 our treatment  does not contain the condition of
vanishing the $GL_{q}(2)$-invariant
\be
     C_2 \equiv tr_{q} (\Omega^2) = 0,
\ee
which is not gauge covariant. Here we have used the
 notion of the $q$-deformed covariant
trace  [4, 8].

 One can rewrite the commutation relations for the gauge
 covariant differential algebra
in terms of the $R$-matrix. Direct checking leads to
the next formulae
\bea
\ba{ll}
 R_{12} dT_{1} T_{2} = T_{2} dT_{1}  R_{12},    \label {f4algdif} &     \\
 R_{12} \Omega_{2} R^{-1}_{12} \Omega _1 + \dfrac{1}{q^2}
    \Omega_1 R_{12} \Omega_2 R^{-1}_{12}      &      \\
  -    \dfrac {q}{1+q^2 +q^4} (E_{12} + (q+\qin) {\cal E}_{21}) tr_q
                \Omega^2 = 0 , &
\ea
\eea
where
\bea
\ba{lr}
E^{ij}_{kl}=\delta^i_k\delta^j_l, &
                      { \cal E}^{ij}_{kl} = \vareps^{ij} \vareps_{kl}.
\ea
 \eea

To construct an exterior differential it is convenient to use the
 definition  based on the dualism between the
exterior algebra of differential forms and the $q$-Lie algebra
of vector fields. In this way
the Leibnitz rule is followed straightforwardly and
it depends only on a special choice  of the $q$-Lie algebra.

Let us start from a general 3D $q$-Lie
 algebra of left-invariant vector fields
    $D^a = (D^{++}, D^{--}, D^{0})$ on the quantum group \su with
    a Lie bracket
\bea
[D^a , D^b]_{B} \equiv D^a D^b - B^{abcd} D^c D^d   = C^{abc} D^c .
\eea
We consider the braiding matrix $ B^{abcd} $ to be unitary,
so that it generates a representation of the permutation group.
Thus, one can easily define the alteration rules for the
tensor algebra of vector fields. Moreover, a generalized Jacobi
identity will be available as well.

The basic left-invariant differential 1-forms $\omega^a$ are defined
 as dual
 objects by means of the scalar product
\bea
\omega^a (D^b) = \delta ^{ab}\;.
\eea
The action of the exterior differential on arbitrary functions $f$ and
differential 1-forms $u$ is defined in analogy with the
 classical case \cite {Kob1}
\bea
\ba{lr}
df(D^a) = D^a f, &  \\
du(D^a, D^b) = -\frac{1}{2} (D^a u(D^b) - B^{abcd} D^c u (D^d) -
       u([D^a, D^b]_{B}),& \label {f5rules}    \\
                         du(D^a, D^b) = - B^{abcd} d u(D^c, D^d).&
\ea
\eea
Rules  for the exterior differentiation
 of the differential ($n>1$)-forms can be
generalized in a similar fashion.
 The Cartan-Maurer equations have a standard
form
\bea
d\omega ^d (D^a, D^b) = \frac{1}{2} C^{abc} \omega^d (D^c). \label {cm}
\eea
Note, that although the braiding matrix $ B^{abcd} $ determines a
corresponding exterior $B$-algebra, nevertheless, commutation relations
for the differential 1-forms $\omega^a$ are not specified completely.

    As a concrete example we consider
 the $q$-Lie algebra
(\ref{Dalg}) which is consistent with the gauge covariant algebra of
left-invariant differential 1-forms (\ref{f4algdif}). In this case
differentiation rules (\ref{f5rules}) can be rewritten
 in a more familiar form after using the
 explicit tensor representation for the exterior products
of 2-forms $\omega^a \wedge  \omega^b$
\bea
\ba{lr}
   d=\omega^a D^a,   &      \\
   d(\op f)= d\op f +\beta \oo \op D^0 f -{\omega^{(2)0}} D^{++} f,   &  \\
   d(\om f) = d\om f +\beta q^2 \oo \om D^0 f + q^2 {\omega^{(2)0}} D^{--}f,
                         \label {f5rules2}       &  \\
   d(\oo f)=\omega^{(2)0} f+
 \beta q^2 \op \oo D^{--}f +\beta \om \oo D^{++},&\\
   \omega^{(2)0} \equiv \dfrac{1}{1+q^2}(\op\om-q^2\om\op)\,,& \\
 \beta \equiv \dfrac{1+q^4}{q^2(1+q^2)}\,\,.    &
\ea
\eea
   It should be noted that the formulae (\ref{f5rules2}) involve just three
   independent basis differential 2-forms
 $\oo \op,  \oo \om , \omega^{(2)0}$
    in the space of exterior 2-forms as in the classical case.
The fourth linearly independent basis form $\sigma^0$ is defined as
follows
\bea
 \sigma^0 =  \dfrac{1}{1+q^2}(\op\om+q^2\om\op).
\eea
The form $\sigma^0$ takes a non-zero value only for the
symmetrical tensor product $D^0\bigotimes D^0$:
\bea
\sigma^0(D^0,D^0) = \rho,
\eea
where the number $\rho$ vanishes in the classical
 limit $q\rightarrow1$. Due to this
property the form $\sigma^0$ does not appear in eqs. (\ref{f5rules2}).

   Starting from the most general form for the exterior differentiation
of 2-forms and taking into account the  condition
$d^2 = 0$ one can derive the next relations:
 \bea
     \ba{ll}
    d(\op \oo f)=-\dfrac{1}{q^2} {\upsilon} D^{++}f, & \\
    d(\om \oo f)=q^2  {\upsilon} D^{--}f, &\\
    d\omega^{(n>2)}=0, &                     \\
                        d(\omega^{(2)0} f)= \beta  {\upsilon} D^0 f,&\\
                        d(\sigma^0\,f)= \rho {\upsilon} f,& \\
  \upsilon \equiv \dfrac{1}{2}       &
(\oo\op\om-q^2\oo\om\op)\,,
     \ea
 \eea
here, $ \upsilon$  is the volume 3-form.

    Having carried out some calculations one can find the
explicit expressions for
    the Cartan-Maurer equations (\ref {cm})
    \bea
    d\Omega=\Omega^2-\frac{q^2}{1+q^2}{ E}\, tr_q \Omega^2.
           \eea
One can see immediately that the r.h.s. of this equation
contains only traceless part of the $\Omega^2$.

\section{Conclusion}
 \indent

One can try to formulate a gauge Yang-Mills theory for the quantum group
\su in analogy with the covariant $GL_q(N)$ version proposed
in ref. \cite {Isa2}. Differentiation over  the matter scalar
fields $\phi^i$ and gauge field $A^i_j$
 in a case of \su
will be more complicated due to the non-standard Leibnitz rule.
 For instance, one has the following differentiation
 rules for the products of two scalar fields
 \bea
\ba{ll}
 d(\phi^i \phi^j) = d\phi^i \phi^j +\frac{1}{q^2} \phi^i  d\phi^j +
       \frac{q^2-1}{q^3} \vareps^{ij} d\phi_{k} \phi^{k}, &      \\
 d(\bar {\phi}^i \bar {\phi}^j) = d\pbar^i \pbar^j +q^2 \pbar^i d\pbar^j,
                                    &   \\
 d(\pbar^i \phi^j) = d\pbar^i \phi^j + q^2 \pbar^i d\phi^j +
  (q^2-1)\pbar^i \pbar^j d\phi_k\phi^k.   &
\ea
  \eea
  The gauge field $A^i_j$ satisfies the same commutation relations
 that right invariant differential
 1-forms $Z = -\dfrac{1}{q^2} T \Omega T^{-1}$ do.
   Another possible way toward a consistent gauge Yang-Mills theory
   corresponds to the differential calculus with a $q$-Lie algebra
   differed from one defined by eqs. (\ref{Dalg}).
   Search for a suitable $q$-Lie algebra is not carried out
   directly due to possible non-lexicographic differentiation rules.

\end{document}